\documentclass[13pt, reqno]{amsart}
\usepackage{amsmath, amsthm, amscd, amsfonts, amssymb, graphicx, color}

\newtheorem{theorem}{Theorem}[section]

\newtheorem{proposition}[theorem]{Proposition}
\newtheorem{corollary}[theorem]{Corollary}
\theoremstyle{definition}
\newtheorem{definition}[theorem]{Definition}

\theoremstyle{remark}
\newtheorem{remark}[theorem]{Remark}
\numberwithin{equation}{section}

\begin{document}
	%\setcounter{page}{1}
	
	%\centerline{}
	
	%\centerline{}
	
	\title[Kernels of Perturbed Hankel Operators]{Kernels of Perturbed Hankel Operators}
	
	\author[Arup Chattopadhyay and Supratim Jana]{Arup Chattopadhyay$^*$ and Supratim Jana$^{**}$}
	\maketitle

	%\subjclass[2020]{ Mathematics }
	
	\paragraph{\textbf{Abstract}}
	In the classical Hardy space $H^2(\mathbb{D})$, it is well-known that the kernel of the Hankel operator is invariant under the action of shift operator S and sometimes nearly invariant under the action of backward shift operator $S^{*}$. It appears in this paper that kernels of finite rank perturbations of Hankel operators are almost shift invariant as well as nearly $S^*$- invariant with finite defect. This allows us to obtain a structure of the kernel in several important cases by applying a recent theorem due to Chalendar, Gallardo, and Partington.
	
	\vspace{0.5cm}
	\vspace{0.5cm}

	\section{\textbf{Introduction and preliminaries}}
	
	The theory of Toeplitz and Hankel operators has a long and intriguing history and is distinguished by amusing results and deep connections with many fields of mathematics, physics, statistical mechanics,
	and so on (see, for example, \cite{BS,VP}). Despite the fact that Toeplitz operators and Hankel operators are quite different in their nature, they are closely related to each other. It is well known that the kernel of a 
	Toeplitz operator is nearly invariant under the backward shift operator $S^*$ acting on the classical Hardy space $H^2(\mathbb{D})$, and the concept of nearly backward shift invariant subspaces was first
	introduced by Hitt in \cite{HITT, SAR!} as a generalization to Hayashi's results concerning Toeplitz kernels in \cite{EH}, and studied further by Sarason in \cite{SAR}.  Recently, there has been an extensive study concerning nearly invariant subspaces with finite defect under the backward shift operator in both the scalar-valued and vector-valued Hardy
	spaces (see, for example, \cite{CP, CGP, CGP1, CD, CDP1, CDP, ROL}). In this direction, Liang and Partington recently provided a connection between kernels of finite-rank perturbations of Toeplitz operators and nearly invariant subspaces with finite defect under the backward shift operator acting on $H^2(\mathbb{D})$. In other words, they have examined the following question, which is closely related to the invariant subspace problem:
	\vspace{-0.2in}
	
	\begin{equation}\label{eq1}
		\begin{split}
			&	\textup{\emph{Given a Toeplitz operator $T$ acting on the Hardy space $H^2(\mathbb{D})$, is the kernel of a finite-rank}}\\
			& 	\textup{\emph{perturbation of T nearly backward shift invariant with finite defect?}}
		\end{split}
	\end{equation}
	\vspace{-0.2in}
	
	Motivated by the work of Liang and Partington \cite{LP}, we continue to study the above problem \eqref{eq1} in the context of the Hankel operator. More precisely, the purpose of this article is to study the kernels of finite-rank perturbations of Hankel operators in $H^2(\mathbb{D})$ and find its connection with almost shift-invariant subspaces as well as nearly $S^*$- invariant subspaces with finite defect. While investigating the above problem \eqref{eq1} in the setting of the Hankel operator, it turns out that the kernel of a finite rank perturbation of the Hankel operator $H_{\phi}$ is nearly $S^*$-invariant subspaces with finite defect as well as almost shift-invariant subspaces for several important classes of symbols $\phi$ (see, Theorem~\ref{mainthm1} and Theorem~\ref{mainthm3}). In this direction, the readers are also referred to \cite{GU1,GU2} for other properties of Hankel Kernels.
	\vspace{0.01in}

	We now introduce some necessary definitions and notations, after which we shall be able to summarize the paper's main contributions.

	Let $Hol(\mathbb{D})$ denote the space of all holomorphic functions on the open unit disk $\mathbb{D}$. The Hardy Space $H^p (\mathbb{D})$ is defined by  $$H^p(\mathbb{D}) := \{ f\in Hol(\mathbb{D}):\sup_{0\le r<1}\frac{1}{2\pi} \int_0^{2\pi}\lvert f(re^{i\theta})\rvert^{p} d\theta < \infty\}. $$ 
	\par The special case $p=2$, the space is known as the Hardy-Hilbert space $H^2(\mathbb{D})$ and has an equivalent definition  $$H^2(\mathbb{D}):= \{f\in Hol(\mathbb{D}) : f(z)= \sum_{n=0}^\infty a_nz^n \text{ with } ||f||^2 = \sum_{n=0}^\infty |a_n|^2 <\infty \}.$$
	
	It is important to note that both the radial and the non-tangential limit of a function of $H^2(\mathbb{D})$ exist and coincide with each other. A classical result due to Fatou (see \cite{Du}, for instance) states that the radial limit $bf(e^{it}):=\lim_{r\rightarrow1-}f(re^{it})$ exists a.e. on the boundary of $\mathbb{D}$ (that is, on the unit circle $\mathbb{T}$). In this regard, it is well known that the space $H^2(\mathbb{T})$, generated by $\{ e^{int}: n \in \mathbb{N} \text{ and } n\ge0 \}$, is isometrically embedded as a closed subspace of $L^2(\mathbb{T})$ via $$ \sum_{n=0}^{\infty} a_nz^n \rightarrow \sum_{n=0}^{\infty} a_ne^{int},$$ and the relation $||f||_{H^2(\mathbb{D})}=||bf||_{H^2(\mathbb{T})}= ||bf||_{L^2(\mathbb{T})}$ holds ($bf$ also known as the boundary function of $f$ ).
	This gives us a natural orthogonal decomposition of $L^2(\mathbb{T})$ as $L^2(\mathbb{T}) = H^2(\mathbb{T}) \oplus (H^2(\mathbb{T}))^\perp$. The space $(H^2(\mathbb{T}))^\perp$ is denoted by $\overline{{H_0^2}}$ which is generated by $\{ e^{int}: n \in \mathbb{N} \text{ and } n<0 \}.$ 
	
	Let $L^\infty := L^\infty(\mathbb{T})$ be the space of all essentially bounded functions on $\mathbb{T}$ and $H^\infty := H^\infty(\mathbb{D})$ is the Banach algebra of bounded holomorphic functions on $\mathbb{D}$ with the norm defined as $||f||_\infty = sup_{z\in\mathbb{D}} |f(z)|$. It is quite obvious that $ H^\infty (\mathbb{D}) \subseteq H^2(\mathbb{D}) .$ So, in a similar way, the radial (or, non-tangential) limit of $H^\infty$ functions belong to $L^\infty$, and hence $H^\infty$ can be viewed as a Banach subalgebra of $L^\infty.$ A function of $H^\infty$ whose boundary value is unimodular on $\mathbb{T}$ is called an inner function. A function $f\in H^1(\mathbb{D})$ with  the following expression  $$f(re^{i\theta}) = \alpha \exp\left({\frac{1}{{2\pi}}} \int_0^{2\pi} \frac{e^{it} + re^{i\theta} }{e^{it} - re^{i\theta} } ~k(e^{it})dt \right), \quad \text{for}~ re^{i\theta}\in \mathbb{D},  $$ is called an outer function, where k is a real-valued integrable function defined on $\mathbb{T}$ and $\alpha$ is a unimodular constant. It is very well known that by inner-outer factorisation theorem any non-zero function $f$ of $ H^1(\mathbb{D})$ has a representation $f= I \cdot O $, where $I$ is called an inner part and $O$  is called an outer part, and this factorisation is unique up to a constant of modulus 1 (cf. \cite{CGP1}). 
	
	Let $ P: L^2(\mathbb{T}) \rightarrow H^2(\mathbb{T})$ be the orthogonal projection on $H^2(\mathbb{T})$ defined by $$(Pf)(z) = \int_{\mathbb{T}}{\frac{f(\xi)}{1-\overline{\xi}z}~dm(\xi)},~ |z|<1. $$ Next, we define the flip operator $J$ as follows
	$$J:L^2(\mathbb{T}) \rightarrow L^2(\mathbb{T}) \text{ by } (Jf)(\xi) = f(\overline{\xi}), \text{ for } \xi \in \mathbb{T}.$$ The action on $f$ by the flip operator $J$ will be denoted by $\breve{f} ~(:=Jf)$ in the later parts. Also, from now on, $H^2(\mathbb{D})$ and $H^2(\mathbb{T})$ will be meant by simply $H^2$ as desired in the phrases later. Given an $L^\infty$ function $\phi$, the Toeplitz and Hankel operators are defined respectively as follows $$ T_\phi : H^2 \rightarrow H^2 \text{ by } T_\phi(f)= P(\phi f) \text{ for all } f \in H^2,$$ and 
	$$ H_\phi : H^2 \rightarrow H^2 \text{ by } H_\phi (f) = PJ(\phi f) \text{ for } f \in H^2.$$ For an inner function $\theta$, the model space is denoted by $K_\theta$ and is defined by $K_\theta = H^2 \ominus \theta H^2 = H^2 \cap \theta \overline{H_0^2}.$ It is well known that $||T_\phi|| = ||\phi||_\infty = ||H_\phi||, \text{ and } J^* = J= J^{-1}  .$ The reader is referred to ( \cite{ARO}, \cite{NKN}, \cite{VP}, \cite{BH}, \cite{GRM} ) for detailed investigations on Toeplitz and Hankel operators. Going further, given a closed subspace $\mathcal{M}\subseteq H^2$, we denote by $P_{\mathcal{M}}: H^2\rightarrow \mathcal{M}$ the orthogonal projection of $H^2$ onto $\mathcal{M}$, and the notation will be used frequently in the subsequent sections for various subspaces $\mathcal{M}$.  
	
	The unilateral shift operator $S:H^2 \rightarrow H^2$  is defined by $Sf(z) = zf(z)$, $z\in \mathbb{D}.$ And the backward shift operator $ S^* : H^2 \rightarrow H^2 $ is defined by $S^*f(z) = \frac{f(z)-f(0)}{z} $ for $f\in H^2, z\in \mathbb{D}.$ A well-known result due to Beurling \cite{Bu} states that if $M$ (closed) is a $S$-invariant subspace of $H^2$, then $M$ can be represented as $M=\theta H^2$, where $\theta$ is an inner function, and hence the subspace $(\theta H^2)^\perp$ $( = K_\theta )$ is invariant under the backward shift operator $S^*$. There are other kinds of subspaces that got plenty of attention in recent studies, which are nearly $S^*$-invariant subspaces and almost invariant subspaces for $S$, defined as follows. 
	\begin{definition}
		(i) A closed subspace $M \subseteq H^2$ is called nearly $S^*$-invariant if $S^*f \in M$ whenever $f\in M$ and $f(0)=0 $ . \\(ii) Furthermore, a closed subspace $M \subseteq H^2$ is said to be  nearly $S^*$-invariant subspace with finite defect $m$ if there is an $m$-dimensional subspace $F$ such that $S^*f \in M+F$ whenever $f\in M \text{ with } f(0)=0;$ and we call $F$ the defect space.\\ (iii) And a closed subspace $M \subseteq H^2$ is said to be almost shift invariant if $S(f) \in M+F \text{ whenever } f \in M $ (that is, $SM \subseteq M+F$) for some finite dimensional subspace $F$.
	\end{definition}
	
	One can check that the kernel of the Toeplitz operator is nearly $S^*$-invariant, and hence by Hitt's theorem \cite{HITT} (see also \cite{SAR})) it has a structure of the form $gK_\theta$, where $\theta$ is an inner function and $g$ is the normalized reproducing kernel at $0$ in the kernel itself. Whereas the kernel of the Hankel operator is shift-invariant, and by virtue of Beurling's theorem, it has the form $\theta H^2$ for some inner function $\theta$. Now, if we consider the subspace $M=\theta H^2$, where $\theta$ is an inner function and $\theta(0) \neq 0$, then $f\in M \text{ with } f(0)=0,$ which immediately implies  that $f=z\theta f_1 \text{ for some function } f_1 \in H^2$. Thus $S^*f = \theta f_1 \in M,$ and hence $M $ is a nearly $S^*$-invariant. And if $\theta(0)=0, $ then $\theta H^2 = z^k \theta_1H^2$ for some finite $k\in \mathbb{N}$ and $\theta_1$ is an inner function with $\theta_1(0) \neq 0$. Therefore, $M$ has the form $z^k$ multiple of a nearly $S^*$-invariant subspace, and hence $M$ is nearly $S^*$-invariant with defect 1. Moreover, every nearly $S^*$-invariant subspace is also almost-invariant for $S$ with defect 1, and it was first observed by Chalendar, Gallardo, and Partington in \cite{CGP}. Getting motivated by this, we investigate the following invariant subspace problem in this article.
	\vspace{0.01in}
	
	\par 
	\textit{ Given a Hankel operator $H_\phi$ acting on the Hardy-Hilbert space $H^2$, is the kernel of an n-rank perturbation of $H_\phi$ nearly $S^*$-invariant subspace with finite defect?  And at the very same time is it an almost shift-invariant subspace? }
	\vspace{0.01in}
	
	We recall that a finite $n$-rank operator $T: \mathcal{H} \rightarrow \mathcal{H}$ acting on a Hilbert space $\mathcal{H}$ has the structure $$T(h)= \sum_{i=1}^n \langle h,u_i \rangle v_i , \text{ for all } h\in \mathcal{H},$$ 
	where $\{u_i\} $ and $\{v_i\}$ are orthonormal sets in $\mathcal{H}.$ An $n$-rank perturbation of Hankel operator $H_\phi : H^2 \rightarrow H^2$ is denoted by $R_n: H^2 \rightarrow H^2$ and is defined by \begin{equation}\label{eqar4}
		R_n(h)= H_\phi (h) +T(h
		) = H_\phi (h) + \sum_{i=1}^n \langle h,u_i \rangle v_i 
	\end{equation} with orthonormal sets $\{u_i\}$ and $\{v_i\}$ in $H^2.$

	The remainder of this paper is organized as follows. In section 2, we study the kernel of the perturbed Hankel operator for several wide classes of symbols and got an affirmative answer to the previous question. We also provide the defect spaces for each classes. This enables us to utilize the C-G-P \cite{CGP} theorem to get the structure of $ker R_1$. This has been shown in section 4. The existence of almost shift invariant subspaces is not entirely obvious to see. But however, in section 3, we manage to prove that the kernel of finite rank perturbed Hankel operator attains this property for the same class of symbols as mentioned in section 2.
	
	\section{\textbf{Nearly $S^*$-invariant $Ker R_n$ with finite defect}}
	
	In this section, we demonstrate that the kernel of the operator $R_n$ is nearly $S^*$-invariant with finite defect for plenty of different important classes, especially recognize the finite dimensional defect spaces. Recall that for $z \in \mathbb{T}$, and $\phi \in L^2$, we denote $ \breve{\phi}(z)=J(\phi)(z)= \phi(\overline{z}) $. First, we start with the following theorem which will be needed to prove our main results in this section as well as in later sections.
	\begin{theorem}\cite[Theorem 4.5.4 ]{ARO}
		Let the functions $\phi \text{ and } \psi$ be in $L^\infty$. Then $$T_{\breve{\phi}} H_{e^{i\theta}\psi} + H_{e^{i\theta}\phi} T_\psi = H_{e^{i\theta}\phi \psi} .$$ 
	\end{theorem}
	
	\noindent As a consequence of the above theorem, if $\psi$ is in $H^\infty$, then  we have $$T_{\breve{\psi}}H_\phi = H_{\phi \psi} = H_{\phi}T_{\psi}. $$
	Therefore, in particular, considering $\psi = e^{i\theta} = z$, we get
	\begin{equation} \label{eq6}
		T_{\overline{z}}H_\phi = H_{\phi z}.
	\end{equation} 
	If $h \in Ker R_n$, then the following is satisfied
	\begin{equation}\label{eq7}
		H_\phi (h) + \sum_{i=1}^n \langle h,u_i \rangle v_i =0.
	\end{equation}
	Applying $T_{\overline{z}}$ on the both sides of \eqref{eq7} and  using \eqref{eq6}, we conclude
	\begin{equation*}
		H_{\phi z}(h) + \sum_{i=1}^n \langle h, u_i \rangle S^* v_i=0.
	\end{equation*}
	Now, let $h \in Ker R_n$ be such that $h(0)=0$, then the equation \eqref{eq7} implies the following
	\begin{equation}\label{eq5}
		H_{z \phi}\left(\frac{h}{z}\right) + \sum_{i=1}^n \langle h, u_i \rangle v_i=0
		\Longleftrightarrow
		J\left( \phi z \frac{h}{z}\right) + \sum_{i=1}^n \langle h, u_i \rangle v_i \in \overline{H_0^2}.
	\end{equation}
	
	Therefore, the query of nearly $S^*$-invariant with finite defect is that, \textit{for each $h\in Ker R_n$ with $h(0)=0$, find a vector w in some suitable finite dimensional space F such that } $$S^*(h) +w = \frac{h}{z} + w \in Ker R_n,$$ which is equivalent to the following equation. 
	\begin{equation}\label{eqar1}
		H_{\phi}\left(\frac{h}{z}+w \right) + \sum_{i=1}^n \left\langle (\frac{h}{z} +w), u_i \right\rangle v_i =0.
	\end{equation}
	Again, by applying $T_{\overline{z}}$ on the both sides of \eqref{eqar1}, we obtain
	\begin{equation*}
		H_{z \phi}(\frac{h}{z} +w) + \sum_{i=1}^n \left\langle (\frac{h}{z} +w), u_i \right\rangle S^* v_i =0  
	\end{equation*}
	\begin{equation}\label{eq8}
		\Longleftrightarrow
		J\left({z \phi} \left(\frac{h}{z} +w \right) \right) + \sum_{i=1}^n \left\langle (\frac{h}{z} +w), u_i \right\rangle S^* v_i \in \overline{H_0^2}.
	\end{equation}
	Now, we try to detect the defect space for some wide important classes of symbols.
	\subsection{$ \phi=0$ a.e. on $\mathbb{T}$}
	
	In this case, $R_n$ is an $n$-rank operator on the Hardy-Hilbert space, and the equation \eqref{eq7} with $\phi$ =0 implies $$Ker R_n = \bigcap_{i=1}^n\left( \bigvee \{u_i\} \right)^ \perp = H^2 \ominus \left( \bigvee \{u_i, i=1,2,...,n\} \right), $$
	where $\bigvee$ denotes the closed linear span in $H^2.$ So, whenever $h \in Ker R_n$ with $h(0)=0$, we always have $$S^*h \in KerR_n \oplus  \left( \bigvee \{u_i, i=1,2,...,n\} \right) = H^2,$$ which leads us to the following observation on the nearly $S^*$-invariant subspace with finite defect.
	\begin{proposition}
		Let $\phi = 0 $ a.e. on $\mathbb{T}.$ Then the subspace $Ker R_n$ is a nearly $S^*$-invariant subspace with finite defect, the dimension of the defect space is at most $n$, and the defect space is contained in the subspace $ F = \bigvee \{u_i, i=1,2,...,n\}.$
	\end{proposition}
	
	\begin{remark}
		When we consider the symbol $\phi$ to be in $H_0^2 ~(=zH^2),$ we get the same defect space as above since  in that case also $ H_\phi (h) = 0 $ for all $h \in H^2.$
	\end{remark}

	\subsection{ $\phi = \overline{\theta} g$, where $ g \in \mathcal{G}( H^\infty) $ }
	
	The celebrated Douglas-Rudin theorem states that $L^\infty = clos_{ L^\infty }\{\overline{\Theta}H^\infty:~~\Theta\in \Omega\}$, where $\Omega$ is the collection of all inner functions (for more details see \cite{BS09, JB, KMD,  NKN}). In other words, the collection $\{\sum\limits_{i=1}^n \bar{\xi_i}\eta_i,~~i=1,2,\ldots,n; ~\xi_i,\eta_i\in \Omega\}$ is  dense  in $L^\infty.$  In this section, we identify the defect space corresponding to the symbol $\phi = \overline{\theta} g$, where $ g \in \mathcal{G}(H^\infty):=$ $\{g\in H^{\infty}:~g^{-1}( =\tilde{g}, \text{ say } )  \in H^\infty,~|g|=\text{constant}\}$ , a subclass of $(\overline{\Theta}H^\infty).$ At the end, as a consequence of our main theorem we obtain a few corollaries related to some special classes of symbols.

The relation \eqref{eq5} gives 
	\begin{equation}\label{eq12}
		\text{ ( say ) } \overline{f_1} = J\left( \overline{\theta} g z \frac{h}{z}\right) + \sum_{i=1}^n \langle h, u_i \rangle v_i \in \overline{H_0^2},
	\end{equation}
	which by applying $J$ on the both sides of \eqref{eq12} implies $ \frac{h}{z}+ \sum\limits_{i=1}^n \langle h, u_i \rangle \theta \tilde{g} J({z v_i}) \in \theta \tilde{g}H^2, $ and hence $ \sum\limits_{i=1}^n \langle h, u_i \rangle \theta \tilde{g} J({z v_i}) \in H^2 .$ On the other hand the equation \eqref{eq8} yields the following 
	
	$$J\left({z \overline{\theta} g} \left(\frac{h}{z} +w \right) \right) + \sum_{i=1}^n \left\langle (\frac{h}{z} +w), u_i \right\rangle S^* v_i \in \overline{H_0^2}  .$$

	\begin{equation} \label{eq13}
		\Longleftrightarrow J\left({z \overline{\theta} g} \left(\frac{h}{z} +w \right) \right) + \sum_{i=1}^n  \frac{1}{\mid g \mid ^2}  \left\langle J\left( z\overline{\theta} g (\frac{h}{z} +w) \right),J\left( z\overline{\theta} g u_i \right) \right\rangle S^* v_i \in \overline{H_0^2}.
	\end{equation}
	
	Now we decompose $f_{1}$ in the following manner. 
	
	Suppose $g u_i = y_i$, then the decomposition for $y_i$ with respect to $\theta$ would be $y_i= y_{i1} + \theta y_{i2}$, where $y_{i1} \in K_{\theta} \text{ and } y_{i2} \in H^2. $ Thus, $J(z \overline{\theta} g u_i) =J( z \overline{\theta} y_{i1}) +J( zy_{i2}),$ where $J(z\overline{\theta} y_{i1}) \perp \overline{H_0^2}$. As $\overline{f_1} \in \overline{H_0^2}$, then  $ \overline{f_1} = \overline{f_{11}} + \overline{f_{12}}$, where $\overline{f_{11}} \in  \bigvee_{i=1}^n \{ J(z y_{i2}) \} (:= J(B),$ where $B:=\bigvee_{i=1}^n \{ zy_{i2} \}),$ and $\overline{f_{12}} \in \overline{H_0^2} \ominus J(B). $    Therefore $\langle \overline{f_{12}}, J(z\overline{\theta}  y_{i1}) \rangle =0 $, and hence $ \langle \overline{f_{12}} , J(z\overline{\theta} g u_{i} )\rangle =0 \text{ for all } i=1,2,...,n .$

	Using the above decomposition we define $w$ as,
	
	\begin{equation}\label{eqar3}
		  \hspace{0.5in} w = \sum_{k=1}^n \langle h,u_k \rangle \theta \Tilde{g} J ({z v_k}) - \overline{z} \theta \Tilde{g} J(\overline{f_{11}}) \in H^2,
	\end{equation}

	(since  $\overline{z} \theta \Tilde{g} J(\overline{f_{11}}) = \tilde{g} P_{\theta H^2} (gu_i) \in H^2 $, and $ \sum\limits_{i=1}^n \langle h, u_i \rangle \theta \tilde{g} J({z v_i}) \in H^2. $)
	
	 $$ \Longleftrightarrow   w = \sum_{k=1}^n \langle h,u_k \rangle H_{ \breve{ \theta}\breve{ \Tilde{g} } } ( z v_k ) - \Tilde{g} \overline{z} ( \theta J(\overline{f_{11}}) ). $$ 
	
Now, substituting $w$, achieved in \eqref{eqar3},  in $J\left({z \overline{\theta} g} \left(\frac{h}{z} +w \right) \right)$, we have 
	\begin{align}
		\nonumber J\left({z \overline{\theta} g} \left(\frac{h}{z} +w \right) \right)=J \left( \overline{\theta}g z \frac{h}{z} \right) +  \sum_{i=1}^n \langle h,u_i \rangle v_i - \overline{f_{11}} \\
		\nonumber = J \left( \overline{\theta}g z \frac{h}{z} \right) +  \sum_{i=1}^n \langle h,u_i \rangle v_i - \overline{f_1} + \overline{f_{12}} = \overline{f_{12}} \text{ ( using \eqref{eq12} )},
	\end{align}
which further implies
	 \begin{align} 
		\nonumber J\left({z \overline{\theta} g} \left(\frac{h}{z} +w \right) \right) + \sum_{i=1}^n  \frac{1}{\mid g \mid ^2}  \left\langle J\left( z\overline{\theta} g (\frac{h}{z} +w) \right),J\left( z\overline{\theta} g u_i \right) \right\rangle S^* v_i \\  
		\nonumber = \overline{f_{12}} + \sum_{i=1}^n \frac{1}{|g|^2}\langle \overline{f_{12}}, J ( z g \overline{\theta}u_i ) \rangle S^*v_i = \overline{f_{12}} \in \overline{H_0^2}. 
	\end{align} 
	Hence, \eqref{eq13} is satisfied. Therefore, in conclusion, the above discussions leads us to the  following theorem.
	
	\begin{theorem}\label{mainthm1}
		\textit{  Suppose $\phi = \overline{\theta}g$ for some non-constant inner function $\theta \text{ and } g\in \mathcal{G}(H^\infty)$.  Then the subspace $Ker R_n$ is nearly $S^*$-invariant with defect at most $2n$, and the defect space is contained in the space  $$ F= \bigvee \{ H_{ \breve{ \theta} \breve{ \Tilde{g} } } ( z v_i ) ,~\Tilde{g} ( P_{ \theta H^2 } (g u_i) ):~  \text{ for } i=1,2,\ldots,n \}. $$ }
	\end{theorem}
	
From the above theorem, we have the following two corollaries; first one is replacing $\theta$ by $1$, and the second one is replacing $g$ by $1.$
	
	\begin{corollary}
		
		Suppose $\phi = g$, such that $g \in \mathcal{G}(H^\infty)$. Then  the subspace $Ker R_n$ is nearly $S^*$-invariant with defect at most $2n$, and the largest defect space is 
		$$ F= \bigvee \{ H_{ \breve{ \Tilde{g} } } ( z v_i ),  u_i :~ \text{ for } i=1,2,\ldots,n \}. $$
		
	\end{corollary}

	\begin{corollary} \label{th1}
		Suppose $\phi = \overline{\theta}$ for some non-constant inner function $\theta$ . Then the subspace $Ker R_n$ is  nearly $S^*$-invariant with defect at most $2n$, and the defect space is contained in the space  $$ F= \bigvee \{ H_{ \breve{\theta} }(z v_i) ,  P_{\theta H^2}(u_i):~  \text{ for } i=1,2,\ldots,n \}. $$
	\end{corollary}

	\section{\textbf{{Almost Shift Invariant $ker R_n$ with finite defect}}}
	
	In this section, we prove that the kernel of the operator $R_n$ (as defined in \eqref{eqar4}) is not only nearly $S^*$-invariant subspace with finite defect but also almost shift invariant subspace for the same classes of symbols mentioned in the previous section. In this case, we would also like to identify the finite-dimensional defect space using a process similar to the one used in the last section. To start with, let's recall the definition and a part of the theorem of the previous section.
	
	\begin{definition}
		A closed subspace $M \subseteq H^2$ is said to be almost shift invariant if $Sf \in M+F \text{ whenever } f \in M $ (that is, $SM \subseteq M+F$) for some finite dimensional subspace $F$.
	\end{definition}
	In such a case, the smallest possible dimension of such $F$ is called the defect of the space $M$.
	
	The relation we shall use here is $$ T_{\overline{z}}H_\phi = H_{\phi z}, \text{ for } z \in \mathbb{T} ~~\text{(see equation \eqref{eq6}}). $$ 
	If $h \in Ker R_n$, then from \eqref{eq7} we have $H_\phi (h) + \sum\limits_{i=1}^n \langle h,u_i \rangle v_i =0$, and this immediately implies $ H_{\phi z}(h) + \sum\limits_{i=1}^n \langle h, u_i \rangle S^* v_i=0 $ by applying $T_{ \overline{z} }$ on the both sides of the last equation. We can now rewrite this as follows
	\begin{equation*}
		H_{\phi }(zh) + \sum_{i=1}^n \langle h, u_i \rangle S^* v_i=0.
	\end{equation*}
	\begin{equation}\label{eq15}
		\Longleftrightarrow J( \phi zh ) + \sum_{i=1}^n \langle h,u_i \rangle S^* v_i \in \overline{ H_0^2. }
	\end{equation}
	Therefore, the question of almost shift invariant subspace turns down to the following, \\ \textit{for each $h \in Ker R_n $, find a vector $w$ in an appropriate finite dimensional space $F$ such that $$ Sh +w = zh + w \in Ker R_n, $$} 
	which implies the following equation 
	\begin{equation*}
		H_\phi ( zh + w ) + \sum_{i=1}^n \langle zh+ w , u_i \rangle v_i = 0.
	\end{equation*}
	\begin{equation}\label{eq16}
		\Longleftrightarrow J (\phi( zh+ w )) + \sum_{i=1}^n \langle zh+w, u_i \rangle v_i \in \overline{H_0^2.}
	\end{equation}
	
	Now, we try to detect the defect space for the same classes of symbols used in the last section. 
	
	\subsection{$ \phi=0$ a.e. on $\mathbb{T}$}
	
	In this case, $R_n$ is an $n$-rank operator and the equation \eqref{eq7} with $\phi$ =0 implies $$Ker R_n = \bigcap_{i=1}^n\left( \bigvee \{u_i\} \right)^ \perp = H^2 \ominus \left( \bigvee \{u_i, i=1,2,\ldots,n\} \right).$$
	 So, whenever $h \in Ker R_n$ , we always have $$ Sh \in KerR_n \oplus  \left( \bigvee \{u_i, i=1,2,\ldots,n\} \right) = H^2,$$ which leads us to the following observation concerning the almost shift invariant subspace.
	
	\begin{proposition}
		Let $\phi = 0 $ a.e. on $\mathbb{T}.$ Then the subspace $Ker R_n$ is an almost shift invariant subspace, and the dimension of the defect space is at most $n$ with the defect space being a subspace of the space $ F = \bigvee \{u_i, i=1,2,\ldots,n\}.$
	\end{proposition}

\begin{remark}
	The above result holds same for the symbol $\phi \in zH^2.$
\end{remark}
	
	\subsection{$ \phi= \overline{\theta} $, where $\theta$ is a non constant  inner functions }
	
	The relation \eqref{eq15} gives \begin{equation}\label{eq17}
		\text{ ( say ) } \overline{f_1} = J\left( \overline{\theta}  zh\right) + \sum_{i=1}^n \langle h, u_i \rangle S^* v_i \in \overline{H_0^2},
	\end{equation} 
which by applying $J$ on the both sides of \eqref{eq17} implies that $\sum\limits_{i=1}^n \langle h, u_i \rangle \theta J(S^* v_i) \in H^2.$
On the other hand the equation \eqref{eq16} yields the following 
	
	$$J\left({ \overline{\theta}} \left(z h +w \right) \right) + \sum_{i=1}^n \left\langle (zh +w), u_i \right\rangle  v_i \in \overline{H_0^2}$$ .
	
	\begin{equation}\label{eq18}
		\Longleftrightarrow J\left({ \overline{\theta}} \left(z h +w \right) \right) + \sum_{i=1}^n \left\langle J(\overline{\theta}(zh +w)),J(\overline{\theta} u_i) \right\rangle  v_i \in \overline{H_0^2}.
	\end{equation}

	Next, let $f_1$ be decomposed in the following way.
	
	The decomposition for $u_i$ with respect to $\theta$ would be $u_i= u_{i1} + \theta u_{i2}$, where $u_{i1} \in K_{\theta}$ and $u_{i2} \in H^2. $ Therefore, $ J(\overline{\theta}u_i) =  J(\overline{\theta} u_{i1}) + J(u_{i2}),$ where $J(\overline{\theta} u_{i1}) \perp \overline{ {H_0^2}}. $ Now, $\overline{f_1} =\overline{ f_{11}} +\overline{f_{12}},$ where $ \overline{f_{11}} \in  \bigvee_{i=1}^n \{ J(u_{i2}-\langle u_{i2}, 1 \rangle) \} (:=J(B), \text{ say} ),~ \overline{f_{12}} \in \overline{H_0^2} \ominus J(B), ~\text{ and }  B:=\bigvee_{i=1}^n \{ u_{i2}-\langle u_{i2}, 1 \rangle \} . $   Therefore $\langle \overline{f_{12}},J( \overline{\theta}  u_{i1} )\rangle =0 ,$ and hence $\langle \overline{f_{12}} , J(\overline{\theta} u_{i}) \rangle =0 \text{ for all } i=1,2,\ldots,n .$ 
	
	Now, using the above decomposition we define $w$ as follows
	
	\begin{align} \label{eq45}
		 w = \sum_{k=1}^n \langle h,u_k \rangle \theta J( S^* v_k ) - \theta J(\overline{f_{11}}) \in H^2,
	\end{align}

	( since $J(\overline{f_{11}}) \in B (\subseteq H^2),$ and $\sum\limits_{k=1}^n \langle h,u_k \rangle \theta J( S^* v_k ) \in H^2.$ )

	$$\Longleftrightarrow w = \sum_{k=1}^n \langle h,u_k \rangle H_{ \breve{ \theta}} ( S^* v_k ) - \theta J(\overline{f_{11}}) . $$
	
Now, substituting $w$, achieved in \eqref{eq45} in $J( \overline{\theta}( z h + w ) )$, we have
	
	$$ J( \overline{\theta}( z h + w ) )=J(\overline{\theta} z h ) + \sum_{i=1}^n \langle h, u_i \rangle S^*v_i - \overline{f_{11}} = J(\overline{\theta} z h ) + \sum_{i=1}^n \langle h, u_i \rangle S^*v_i - \overline{f_1} + \overline{f_{12}} = \overline{f_{12}},$$
	which further implies $$ J\left({ \overline{\theta}} \left(z h +w \right) \right) + \sum_{i=1}^n \left\langle J(\overline{\theta}(zh +w)),J(\overline{\theta} u_i) \right\rangle  v_i = \overline{f_{12}} + \sum_{i=1}^n \langle \overline{f_{12}}, J(\overline{\theta} u_i) \rangle = \overline{f_{12}} \in \overline{H_0^2}.  $$
	
Hence, the equation \eqref{eq18} has been satisfied. Therefore, the above discussions yield to the following conclusion.
	\begin{theorem}\label{mainthm2}
		Suppose $\phi = \overline{\theta}$ for some non-constant inner function $\theta$. Then $Ker R_n$ is an almost shift invariant subspace with dimension of the defect space is at most $2n$, and the defect space is contained in the space $$ \bigvee \{ H_{ \breve{ \theta}} ( S^* v_i ), P_{\theta H^2}(u_i) - \langle \frac{ P_{\theta H^2}(u_i)}{\theta}, 1 \rangle \theta : i=1,2,...,n \}. $$
	\end{theorem}

Next, corresponding to the symbols $\overline{\theta}g$ ($\theta$ is an inner function and $g\in \mathcal{G}(H^{\infty})$), by repeating the same calculations as done before the statement of 
Theorem ~$\ref{mainthm2}$, we reach the following conclusion

	\begin{theorem}\label{mainthm3}
		Suppose  $\phi = \overline{\theta}g $,  where $\theta$ is a non-constant inner function and $g\in H^\infty $ such that $ g^{-1} (=\Tilde{g}, ~\text{say}) \in H^\infty. $  Then $Ker R_n$ is an almost shift invariant subspace with defect at most $2n$, and the defect space is contained in the space
		$$ \bigvee \{ H_{ \breve{\theta} \breve{\Tilde{g}} } ( S^* v_i ), T_{\Tilde{g}}(  P_{ \theta H^2 } (gu_i) - \langle \frac{P_{ \theta H^2 } (gu_i)}{\theta}, 1 \rangle \theta) : i = 1,2,\ldots,n \}. $$
	\end{theorem}

	\section{\textbf{{ Application of C-G-P Theorem }}}
	
	This section contains some application of recent work by Chalender-Gallardo-Partington \cite{CGP} to give a structure of the kernels of rank one perturbed Hankel operators in terms of backward shift-invariant subspaces. We shall consider $n=1$ and denote the operator by $$ R_1 = H_\phi + \langle . , u \rangle v $$ with $||u|| = 1$ and $||v||=1$ such that $S^*v \neq 0.$ First of all, we state the desired theorem (known as C-G-P theorem) for $m$ dimensional defect space.
	\begin{theorem}\cite[Theorem 3.2]{CGP}
		Let M be a closed subspace that is nearly $S^*$-invariant with defect $m$. Then: \\ (1) in the case where there are functions in M that do not vanish at 0, then $$ M = \{ f \in H^2 | f(z)=k_0(z) f_0(z) + z\sum_{i=1}^m k_i(z)e_i(z) : (k_0,\ldots,k_m) \in K \}, $$
		where $f_0$ is the normalized reproducing kernel for M at 0, $${e_1,\ldots,e_m}$$ is an orthonormal basis of F and K is a closed $ S^* \oplus \cdots\oplus S^*$-invariant subspace of the vector-valued Hardy space $H^2(\mathbb{D}, \mathbb{C}^{m+1})$, and $||f||^2 = \sum_{i=0}^m || k_i ||^2. $ 
		\\ (2) In the case where all functions in M vanishes at 0, then $$ \{ 
		f\in H^2| f(z)= z\sum_{i=1}^m k_j(z)e_j(z) : ( k_1,\ldots,k_m ) \in K \} $$ with the same notation as in (1), except K is now closed $ S^* \oplus \cdots\oplus S^*$-invariant subspace of the vector valued Hardy space $H^2(\mathbb{D}, \mathbb{C}^{m})$, and $||f||^2 = \sum_{i=1}^m || k_j ||^2. $ \\
		Conversely, if a closed subspace $M\subseteq H^2$ has a representation as in $(i)$ or $(ii)$, then it is a nearly $S^*$- invariant subspace of defect $m$.
	\end{theorem}
	
	Next, we will be applying the fore stated theorem to represent the kernel of the $R_1$ operator for several symbols and try to find the $S^*$-invariant subspace $K$ as large as possible so that $$ {S^*}^nk_0f_0 + z\sum_{i=1}^m {S^*}^nk_ie_i \in M \text{ or } z\sum_{i=1}^m {S^*}^nk_ie_i \in M $$ for all $n\in \mathbb{N}.$
	
	\subsection{ $\phi=0$ a.e. on $\mathbb{T}$ } In this case $M=Ker R_1 = H^2 \ominus \{ u \}$, a vector hyperplane, which is a nearly $S^*$-invariant subspace with defect $1$ and the defect space is $F=\bigvee \{u\}. $
	For this particular symbol, we would like to refer to \cite{LP}, where the authors already described the structure of $K$, and for reader convenience, we just quote their result as follows.
	\begin{corollary}\cite[Corollary 3.3]{LP}
		Given a nearly $S^*$-invariant subspace $M= H^2\ominus \bigvee \{u\}$ with defect 1, let $f_0 = P_M 1 = 1-\overline{ u(0)}u, v_0 = P( u-u(0){\mid u\mid}^2 ) \text{ and } v_1 = P( \overline{z} {\mid u \mid}^2 )$.  Then  \\  (1) in the case $P_M1 \neq 0$, we have $$ M=\{ f\in H^2 : f=k_0f_0 + k_1zu: ( k_0,k_1 ) \in K \}, $$ with an $S^*\oplus S^*$-invariant subspace $$ K=\{ ( k_0,k_1 ) : \langle k_0, z^n v_0 \rangle + \langle k_1, z^n v_1\rangle =0 \text{ for } n\in \mathbb{N} \}. $$
		(2) In the case $P_M1 =0 $, we have $$ M= \{ f\in H^2 : f= k_1zu : k_1 \in K \} $$ with an $S^*$-invariant subspace $K=\{ k_1 : \langle k_1, z^nv_1 \rangle \text{ for } n\in \mathbb{N} \}.$
	\end{corollary}
	
	\subsection{ $\phi = \overline{\theta} $ with $\theta $ being a non-constant inner function } 
	This is of particular interest since it has a link to shift-invariant as well as nearly $S^*$-invariant subspaces. Now, if  $h \in Ker R_1$, then we have $$H_{\overline{\theta}} (h) + \langle h,u \rangle v =0$$
	$$ \implies J(\overline{\theta}h) + \langle h,u \rangle v \in \overline{H_0^2},$$
	$$\implies h +\langle h,u \rangle \theta Jv \in \theta H_0^2 = z\theta H^2  \implies \theta Jv \in H^2.$$
	On the other hand, we also have $\theta Jv \in \theta \overline{H^2}$, and hence  $\theta Jv \in H^2 \cap z \theta \overline{H_0^2} = K_{z \theta}$ .  Therefore $$M= Ker R_1 \subset z\theta H^2 \oplus \bigvee \{ \theta Jv \}. $$ Consider any vector $h=h_1 + \lambda \theta Jv \in M $ with $h_1 \in z\theta H^2$ and $\lambda \in \mathbb{C}$ such that $ R_1 h =0, $ then it is equivalent to 
	\begin{equation}\label{eq20}
		\lambda ( 1 + \langle \theta Jv, u \rangle ) = - \langle h_1, u \rangle .
	\end{equation} 
	Now we have two separate subsections to represent $ M= Ker R_1 $ in terms of $S^*$-invariant subspaces.
	\subsection{$u \perp z\theta H^2$}  In this case the equation \eqref{eq20} will turn into $  \lambda ( 1 + \langle \theta Jv, u \rangle ) = 0. $ So, either one of them is zero. If $\lambda =0 $, then $ M= z\theta H^2 $ which is clearly a nearly $S^*$-invariant subspace with defect $1.$ Therefore by applying C-G-P theorem we conclude that $$ M=\{ f \in H^2 : f= k_1 z u, \text{ where } k_1 \in H^2 \text{ and } u= \theta \}. $$ Therefore, the defect space is $\bigvee \{ \theta\}$. On the other hand, Corollary~\ref{th1} says that the defect space is a subspace of $ \bigvee \{ H_{\breve{\theta}}(zv),  P_{\theta H^2}(u) \}. $ Here it is indeed, because $ u= u_1 + \theta u_2 \in K_{z\theta} , $ where $ u_1 \in K_\theta (\subset K_{z\theta}). $ Therefore $ \theta u_2 \perp z\theta H^2$ which assures $u_2 $ to be a constant. Hence the conclusion of the Corollary~\ref{th1} gets satisfied. 
	
	For the case $\lambda \neq 0,$ we have $ M=z\theta H^2 \oplus \bigvee \{ \theta Jv \}. $ And lets assume the dimension of the defect space is $2$ and it is generated by $ \{ H_{\breve{\theta}}(zv),  P_{\theta H^2}(u) \}. $ Therefore by applying C-G-P theorem we deduce the following
	$$ P_M (1) = \overline{{\theta}(0) Jv(0)} \theta Jv= \gamma f_0 ~~(\text{here $\gamma$ is  $||P_M(1)||$, and below it is adjusted with $k_0$}).$$
	Therefore, 
	\begin{align*}
		M= \{ f \in H^2 : f=f_0 k_0  + \alpha k_1 z P_{\theta H^2}(u) + \beta k_2 z^2 H_{\breve{\theta}}(zv) \} = z\theta H^2 \oplus \bigvee \{ \theta Jv \},
	\end{align*}
	where $ \alpha = {||P_{\theta H^2}(u)||}^{-1}$, $\beta = {|| H_{\breve{\theta}}(zv) \} ||}^{-1} $, and $ (k_0,k_1,k_2) \in K, $ and $S^* \oplus S^* \oplus S^*$-invariant subspace of $ H^2( \mathbb{D}, \mathbb{C}^3 ). $ By a similar argument, we say that $ P_{\theta H^2}(u) = \theta. $ Now, comparing both sides, we have 
	$$ M= \{ f \in H^2 : f=\overline{{\theta}(0) Jv(0)} \theta Jv k_0  + \alpha k_1 z \theta + \beta k_2 z^2 H_{\breve{\theta}}(zv) \} = z\theta H^2 \oplus \bigvee \{ \theta Jv \} .$$
	$$ \implies k_0 =\text{constant}, k_1 \in H^2, ~\text{and} ~k_2 =0. $$ 
	Therefore, by the above discussions, we obtain the following.
	\begin{corollary}
		Given a nearly $S^*$-invariant subspace $ M= z\theta H^2 \oplus \bigvee \{ \theta Jv \} $ with defect $1$, can be represented as $$ M= \{ f \in H^2 : f=f_0 k_0  +  k_1 z \theta \} , $$  where $ K= \{ f\in H^2 : f= ( k_0,k_1 ) \text{ such that }k_0 \in \mathbb{C}  $ and $k_1 \in H^2\},$ an $S^*\oplus S^*$-invariant subspace of $H^2(\mathbb{D},\mathbb{C}^2) .$
	\end{corollary}
	
	\subsection{$u \in z\theta H^2$}
	In this case, we split $u$ into $ u= u_1 \oplus u_2 $, where $ u_1 \in K_{z\theta} $ and $u_2 \in z \theta H^2 $. Therefore the relation \eqref{eq20} turns down to $$ \lambda ( 1+ \langle \theta Jv, u_1 \rangle ) = -\langle h_1, u_2 \rangle. $$
	According to the Corollary~\ref{th1}, $Ker R_1$ is nearly $S^*$-invariant with at most two-dimensional defect, and the defect space is formed by either of $H_{\breve{\theta}}(zv),  P_{\theta H^2}(u)$, or both. In this case also $M= KerR_1 \subset N= z\theta H^2 \oplus \bigvee \{ 
	\theta Jv \}.$ Based on the value of $( 1+ \langle \theta Jv, u_1 \rangle ) (= w_\theta, \text{ say })$ there are two results discussed bellow. 
	
	For the first case, if $w_\theta =0,$ then $KerR_1 = N \ominus \vee \{u_2\}.$ 
	Therefore, $P_M(1)= P_N-\overline{u_2(0)}u_2= \overline{\theta(0)Jv(0)}\theta Jv - \overline{u_2(0)}u_2 = \alpha f_0$, where $\alpha = ||P_M||. $
	Thus, $M= \big\{ f\in H^2: f=k_0f_0 + zk_1 H_{\breve{\theta}}(zv) + z^2k_2 P_{\theta H^2}(u) \big\}= \big(z\theta H^2\oplus \vee \{ \theta Jv \}\big) \ominus \vee \{ u_2\} $ (here the norms are been taken care of by an initial multiplication with $k_0,k_1,k_2$). 
	$$ \implies M=\big\{ f \in H^2 : f= \overline{\theta(0)Jv(0)}k_0\theta Jv - \overline{u_2(0)}u_2k_0 + k_1 \theta  Jv + z^2 k_2 \theta u_{12} \big\} $$ 
	$$\hspace{-1.7in}= \big(z\theta H^2\oplus \vee \{ \theta Jv \} \big) \ominus \vee \{ u_2\},$$
	$$ \implies \overline{ \theta(0)Jv(0) }k_0 + k_1 \in \mathbb{C}, (\overline{u_2(0)}k_0 +zk_2) u_2 \in z\theta H^2 \text{ as } zP_{\theta H^2}(u)= z\theta u_{12}= P_{z\theta H^2}(u)= u_2, $$
	$$   \text{ and }  \langle k_0, z^n v_0 \rangle + \langle k_2, z^n v_2 \rangle=0, \text{ where } v_0 = \theta(0) Jv(0) u_2 - u_2(0)|u|^2 \text{ and } v_2= \overline{z}|u_2|^2. $$
	Therefore, from the above discussions, we conclude the following.
	\begin{corollary}
		Given $M= \big(z\theta H^2 \oplus \vee\{\theta Jv\}\big) \ominus \vee\{ u_2\},$ a nearly $S^*$-invariant subspace with two dimensional defect can be represented as $$  M=\big\{ f \in H^2 : f= \overline{\theta(0)Jv(0)}k_0\theta Jv - \overline{u_2(0)}u_2k_0 + k_1 \theta  Jv + z k_2 u_2 \big\} $$ with $ \overline{ \theta(0)Jv(0) }k_0 + k_1 \in \mathbb{C}, ~(\overline{u_2(0)}k_0 +zk_2) u_2 \in z\theta H^2, \text{ and } \langle k_0, z^n v_0 \rangle + \langle k_2, z^n v_2 \rangle=0,$ where $v_0 = \theta(0) Jv(0) u_2 - u_2(0)|u|^2 \text{ and } v_2= \overline{z}|u_2|^2.$
	\end{corollary}
	Now for the second case, if $w_\theta \neq 0$, then  $KerR_n $ turns down to $  N\ominus \vee \{ u_2 +w_\theta \theta Jv \}$. Therefore $P_M = \overline{ \theta(0)Jv(0) }\theta Jv - \overline{\rho_\theta}(u_2+\overline{w_\theta}\theta Jv), \text{ where } \rho_\theta = \overline{u_2(0)}+w_\theta \overline{\theta(0)Jv(0)} .$ Let $P_M = \alpha f_0 $, where $\alpha = || P_M ||.$ Therefore, by considering $$ v_0 = \theta(0)Jv(0)u_2\overline{\theta Jv} + \theta(0)Jv(0) \overline{w_\theta(0)}|Jv|^2 - \overline{\rho_\theta}(|u_2|^2 + |w_\theta Jv|^2), $$ $$  v_1 = P(\overline{w_\theta}|Jv|^2), \text{ and } v_2 = P( \overline{z}|u_2|^2 ), $$ we have the following corollary. 
	\begin{corollary}
		Given $M=\big(z\theta H^2 \oplus \vee \{ \theta Jv \} \big) \ominus \vee \{ u_2 + \overline{w_\theta}\theta Jv \},$ a nearly $S^*$-invariant subspace with two dimensional defect can be represented as $$ M= \big\{ f\in H^2: f= \overline{\theta(0)Jv(0)}\theta Jv k_0 + (k_1 - \overline{w_\theta}\rho_\theta k_0)\theta Jv +(zk_2 - \rho_ \theta k_0)u_2 \big\}$$ 
		such that $(k_0,k_1,k_2)\in K,$  an $S^*\oplus S^* \oplus S^*$-invariant subspace of  $H^2(\mathbb{D}, {\mathbb{C}}^3),$ and  $$K= \{ (k_0,k_1,k_2) : k_i \text{ satisfies } \eqref{eqar6} \text{ for } i=0,1,2 \},$$ 
		\begin{equation}\label{eqar6}
			\overline{\theta(0)Jv(0)}k_0 - \rho_\theta \overline{w_\theta} +k_1 \in \mathbb{C} , ~\rho_\theta (k_0 - zk_2)u_2 \in z\theta H^2, ~and ~\langle k_0,z^n v_0 \rangle  +  \langle k_1, z^nv_1 \rangle + \langle k_2, z^nv_2 \rangle =0.
		\end{equation}
	\end{corollary}
	\begin{remark}
		For other symbols $\phi=\overline{\theta}g$ ($\theta$ is an inner function, and $g\in \mathcal{G}(H^{\infty})$), we are not providing any application of C-G-P theorem since the calculations will be very bizarre in this situation.
	\end{remark}
	
	\subsection{Results for Other Equivalent Definition of Hankel Operator}\label{sec5}
	
	In the literature, there is another definition of $Hankel \text{ } Operator$ (see \cite{JRP,SCP}), namely $$\hat{H}_{\phi}:~H^2\rightarrow \overline{H^2_0} \quad \text{defined by}\quad \hat{H}_\phi (f) = (I-P) M_\phi (f), $$ where $P : L^2 \rightarrow H^2$ is the orthogonal projection, and $(H^2)^{\perp}=\overline{H^2_0}$. Suppose  $\phi\in L^\infty$, then with the above definition, we have the following relation 
	\begin{equation}\label{anohank}
		(I-P)S\hat{H}_\phi = \hat{H}_{z\phi}, 
	\end{equation}
	instead of the relation $ S^* H_\phi = H_{z \phi} $, occurred for the previous definition. Therefore, using the relation \eqref{anohank} and proceeding similarly as in the earlier sections, we obtain the following similar sort of results but with different defect spaces (spanned by different set of vectors compared to the earlier one) for the perturbed Hankel operator $$ R_n = \hat{H}_\phi + \sum_{i=1}^n\langle \cdot,u_i  \rangle v_i ,$$ where $\{u_i\}_{i=1}^n \in H^2$, and $\{v_i\}_{i=1}^n \in \overline{H_0^2}$ are unit vectors.

	\begin{theorem}
		Suppose $\phi = \overline{\theta}g$ for some non-constant inner function $\theta \text{ and } g\in \mathcal{G}(H^\infty)$. Then the subspace $Ker R_n$ is a nearly $S^*$-invariant subspace with defect at most $2n$ dimensional, and the defect space is contained in the space $$ \bigvee\{ S^*(\tilde{g}(P_{\theta H^2}(zgu_i)), P(\overline{z}\theta \tilde{g}v_i): i=1,2,\ldots,n \}. $$
	\end{theorem}

	\begin{theorem}
		Suppose the symbol  $\phi = \overline{\theta} $ for some non-constant inner functions $\theta.$  Then the subspace $Ker R_n$ is an almost shift invariant subspace with dimension of the defect space is at most $2n$, and the defect space is contained in the space $$\bigvee\{  P_{\theta H^2}( u_i ), P(\theta  (I-P) (zv_i )) :i=1,2,\ldots,n \}.$$
	\end{theorem}

	\section*{\textbf{{Concluding Remark}}}
	
	{\it It is always a matter of interest to find out the explicit structure of almost shift-invariant subspaces in $H^2$ and find a connection with perturbed Hankel operators, and we leave it as a subject of future investigation.}
	\vspace{0.2in}

	\noindent {\textbf{{Data Availability Statement}}} Data sharing is not applicable to this article as
	no data sets were generated or analyzed during the current study.
	\vspace{0.1in}

	\noindent {\textbf{{Declarations}}}
	\vspace{0.1in}
	
	\noindent {\textbf{{Conflict of interest}}} The authors declare that they have no conflict of interest.

	\section*{\textbf{{Acknowledgment}}}
	
	The authors are very thankful to the referee for pointing out some errors in the initial draft. The authors are extremely grateful to Prof. Jonathan R. Partington and Prof. Yuxia Liang for their valuable comments and suggestions on this work, which led us to include Subsection~\ref{sec5}.  A. Chattopadhyay is supported by the Core Research Grant (CRG), File No: CRG/2023/004826, by the Science and Engineering Research Board(SERB), Department of Science \& Technology (DST), Government of India. The second author gratefully acknowledges the support provided by IIT Guwahati, Government of India.

	\vspace{0.1in}
	
	\noindent $^*$ [A. Chattopadhyay] Department of Mathematics, Indian Institute of Technology Guwahati, Guwahati, 781039, India.\\
	\textit{Email address:}~ arupchatt@iitg.ac.in, 2003arupchattopadhyay@gmail.com.
	\vspace{0.05in}
	
	\noindent $^{**}$ [S. Jana] Department of Mathematics, Indian Institute of Technology Guwahati, Guwahati, 781039, India.\\
	\textit{Email address:}~ supratimjana@iitg.ac.in, suprjan.math@gmail.com.

\end{document}